\documentclass[12pt,twoside]{article} 
\usepackage{amssymb}
\usepackage{amsmath}

\date{} 

\begin{document} 

\centerline {\Large{\bf SOME FACTORIZATIONS IN THE TWISTED GROUP}} 

\centerline{} 

\centerline{\Large{\bf ALGEBRA OF SYMMETRIC GROUPS}} 

\centerline{} 

\centerline{\bf {Milena So\v{s}i\'{c}}} 

\centerline{} 

\centerline{Department of Mathematics} 

\centerline{University of Rijeka} 

\centerline{Radmile Matej\v{c}i\'{c} 2, Rijeka 51000, CROATIA} 

\centerline{} 

\centerline{e-mail:  msosic@math.uniri.hr} 

\centerline{} 

\newtheorem{Theorem}{\quad Theorem}[section] 

\newtheorem{Definition}[Theorem]{\quad Definition} 

\newtheorem{Corollary}[Theorem]{\quad Corollary} 

\newtheorem{Proposition}[Theorem]{\quad Proposition} 

\newtheorem{Lemma}[Theorem]{\quad Lemma} 

\newtheorem{Example}[Theorem]{\quad Example} 

\newtheorem{Remark}[Theorem]{\quad Remark} 

\begin{abstract} 
In this paper we will give a similar factorization as in \cite{4}, \cite{5}, where the autors Svrtan and Meljanac examined certain matrix factorizations on Fock-like representation of a multiparametric quon algebra on the free associative algebra of noncommuting polynomials equiped  with multiparametric partial derivatives. In order to replace these matrix factorizations (given from the right) by twisted algebra computation, we first consider the natural action of the symmetric group $S_{n}$ on the polynomial ring $R_{n}$ in $n^2$ commuting variables $X_{a\,b}$ and also introduce a twisted group algebra (defined by the action of  $S_{n}$ on  $R_{n}$) which we denote by ${\mathcal{A}(S_{n})}$. Here we consider some factorizations given from the left because they will be more suitable in calculating the constants (= the elements which are annihilated by all multiparametric partial derivatives) in the free algebra of noncommuting polynomials.
\end{abstract}



\section{Introduction}
Following the papers \cite{4}, \cite{5} by Meljanac and Svrtan, where an explicit Fock-like representation of a multiparametric quon algebra on the free associative algebra of noncommuting polynomials equiped  with multiparametric partial derivatives (see also \cite{3}) is constructed, our task here is to replace the `nonobvious' matrix level factorizations by `somewhat' simpler algebraic manipulations in a twisted group algebra ${\mathcal{A}(S_{n})}$.\\ 
More general factorizations in braid group algebra we can find in \cite{2}.\\ 
In order to construct ${\mathcal{A}(S_{n})}$ we first consider the natural action of the symmetric group $S_{n}$ on the polynomial ring $R_{n}$ in $n^2$ commuting variables $X_{a\,b}$ and let ${\mathcal{A}(S_{n})}=R_{n}\rtimes {\mathbb{C}}[S_n]$ be the associated (twisted) group algebra. Further, we give some factorizations of certain canonical elements in ${\mathcal{A}(S_{n})}$ in terms of simpler elements of ${\mathcal{A}(S_{n})}$. Then by representing ${\mathcal{A}(S_{n})}$ on the free unital associative complex algebra $\mathcal{B}$ (= the algebra of noncommuting polynomials) by using multiparametric partial derivatives, we obtain more easily some matrix factorizations. Similarly, we can apply some factorizations in ${\mathcal{A}(S_{n})}$ in the problem of computing constants (i.e the elements which are annihilated by all multiparametric partial derivatives) in the algebra $\mathcal{B}$. This will be elaborated in the fortcoming paper. The explicit formulas for basic constants in the subspaces of $\mathcal{B}$ up to total degree four are given in \cite{8}.

\section{The algebra ${\mathcal{A}(S_{n})}$}
Let $S_{n}$ denote the symmetric group on $n$ letters, i.e $S_{n}$ is the set of all permutations of a set $M=\{1,2,\dots ,n \}$ equiped with a composition as the binary operation on $S_{n}$ (where the permutations are regarded as bijections from $M$ to itself). Note that the groups $S_{n}$, $n\ge 3$ are not abelian.\\

\noindent Let ${X=\{X_{a\,b} \mid 1\le a,b\le n \}}$ be a set of $n^2$ commuting variables $X_{a\,b}$ and let ${R_n:={\mathbb{C}}[X_{a\,b} \mid 1\le a,b\le n ]}$ denote the polynomial ring, i.e the commutative ring of all polynomials in $n^2$ variables $X_{a\,b}$ over the set ${\mathbb{C}}$ (of complex numbers), with ${1\in R_n}$ as a unit element of $R_n$.\\
First, let $S_{n}$ act on the set $X$ as follows
\begin{equation}\label{gX} 
         g{\bf.}X_{a\,b}=X_{g(a)\,g(b)}\, g.
\end{equation}
This action of $S_{n}$ on $X$ induces the action of $S_n$ on $R_{n}$ given by    
\begin{equation}\label{gp} 
         g{\bf.}p(\dots, X_{a\,b},\dots)=p(\dots, X_{g(a)\,g(b)},\dots)\, g
\end{equation}
for every $g\in S_n$ and any $p\in R_{n}$.\\

In what follows we are going to study a kind of twisted group algebra, which we denote by ${\mathcal{A}(S_{n})}$ and call it a twisted group algebra of the symmetric group $S_{n}$ with the coefficients in the polynomial ring $R_{n}$.\\ 

Recall that the usual group algebra ${{\mathbb{C}}[S_n]=\left\{\sum_{\sigma\in S_{n}} c_{\sigma} \sigma \mid c_{\sigma}\in {\mathbb{C}}\right\}}$ of the symmetric group $S_n$ is a free vector space (generated with the set $S_n$), where the multiplication is given by
$$\left( \sum_{\sigma\in S_{n}} c_{\sigma} \sigma \right)\cdot \left( \sum_{\tau\in S_{n}} d_{\tau} \tau \right)= \sum_{\sigma, \tau\in S_{n}} (c_{\sigma}d_{\tau})\, \sigma \tau .$$
Here we have used the simplified notation $\sigma \tau$ for the composition $\sigma \circ \tau$, i.e the product of $\sigma$ and $\tau$ in $S_{n}$.\\
Now we define more general group algebra 
\begin{equation}\label{twga} 
         {\mathcal{A}(S_{n})}:=R_{n}\rtimes {\mathbb{C}}[S_n]
\end{equation} 
a twisted group algebra of the symmetric group $S_{n}$ with coefficients in the polynomial ring $R_{n}$.\\
Here\, $\rtimes$ denotes the semidirect product. The elements of the set ${\mathcal{A}(S_{n})}$ are the linear combinations 
$$\sum_{g_{i}\in S_{n}} p_{i}\,g_{i} \hspace{25pt} \textrm{with}\hspace{7pt} p_{i}\in R_{n}$$
and the multiplication in ${\mathcal{A}(S_{n})}$ is given by 
\begin{equation}\label{mtwga} 
        (p_{1}g_{1}) \cdot ( p_{2}g_{2}):=(p_{1}\cdot (g_{1}{\bf.}p_{2}))\, g_{1}g_{2},
\end{equation}  
where $g_{1}{\bf.}p_{2}$ is defined by (\ref{gp}) and $g_{1}g_{2}$ denotes the product of $g_{1}$ and $g_{2}$ in $S_{n}$.\\
It is easy to see that the algebra ${\mathcal{A}(S_{n})}$ is associative but not commutative.\\

\noindent Let 
$${I(g)=\{(a,b)\mid 1\le a < b\le n, \, g(a) > g(b) \}}$$ 
denote the set of inversions of $g\in S_{n}$.\\ 
Then to every $g\in S_{n}$ we associate a monomial in the ring $R_{n}$ defined by 
\begin{equation}\label{Xmon} 
         X_{g}:=\prod_{(a,b)\in I(g^{-1})} X_{a\,b} \left(=\prod_{a < b, \, g^{-1}(a) > g^{-1}(b)} X_{a\,b} \right),
\end{equation} 
which encodes all inversions of $g^{-1}$ (and of $g$ too).\\

\noindent More generally, for any subset ${A \subseteq \{1,2,\dots, n \}}$ we will use the notation 
\begin{equation}\label{Xset} 
         X_{A}:=\prod_{(a,b)\in A \times A, \, a<b} X_{a\,b}\cdot X_{b\,a}=\prod_{(a,b)\in A \times A, \, a<b} X_{\{a,\,b\}},
\end{equation} 
because
\begin{equation}\label{XsetA} 
         X_{\{a,\,b\}}:=X_{a\,b}\cdot X_{b\,a}.
\end{equation} 
\begin{Definition}\label{gtwist}
To each $g\in S_{n}$ we assign a unique element ${g^*}\in {\mathcal{A}(S_{n})}$ defined by 
\begin{equation}\label{gtilda} 
         g^*:=X_{g}\, g
\end{equation}
with $X_{g}$ defined by $(\ref{Xmon})$.
\end{Definition}
In what follows we will use the elements ${g^*}\in {\mathcal{A}(S_{n})}$ defined by (\ref{gtilda}).

\begin{Theorem}\label{tmmi} 
For every\, $g_{1}^*, g_{2}^*\in {\mathcal{A}(S_{n})}$\, we have
\begin{equation}\label{mid} 
         g_{1}^*\cdot g_{2}^*=X(g_{1}, g_{2})\, (g_{1}g_{2})^*,
\end{equation}
where the multiplication factor is given by
\begin{equation}\label{mfact} 
         X(g_{1}, g_{2})=\prod_{(a,b)\in I(g_{1}^{-1})\backslash I((g_{1}g_{2})^{-1})} X_{\{a,\,b\}} \left(=\prod_{(a,b)\in I(g_{1})\cap I(g_{2}^{-1})} X_{\{g_{1}(a),\,g_{1}(b)\}} \right).
\end{equation}
\end{Theorem}
\textbf{Proof.}\hspace{5pt} By using the notations (\ref{gtilda}) and abbreviating\, $g_{1}g_{2}=g$\, we have
\begin{align*}
g_{1}^*\cdot g_{2}^*&=\left(X_{g_{1}}\, g_{1}\right)\cdot \left(X_{g_{2}}\, g_{2}\right)=(X_{g_{1}}\cdot g_{1}{\bf.}X_{g_{2}})\, g = \left(X_{g_{1}}\cdot \prod_{(c,d)\in I(g_{2}^{-1})} X_{g_{1}(c)\,g_{1}(d)}\right) g \\
&=\left(X_{g_{1}}\cdot \prod_{(g_{1}^{-1}(a),g_{1}^{-1}(b))\in I(g_{2}^{-1})} X_{a\,b}\right) g \\
&=\left(X_{g_{1}}\cdot \prod_{(a,b)\in I(g^{-1})\backslash I(g_{1}^{-1})}X_{a\,b}\, \cdot \prod_{(b,a)\in I(g_{1}^{-1})\backslash I(g^{-1})} X_{a\,b}\right) g \\
&=\left(\prod_{(a,b)\in I(g_{1}^{-1})}X_{a\,b}\, \cdot \prod_{(a,b)\in I(g^{-1})\cap I(g_{1}^{-1})}X_{a\,b}^{-1}\, \cdot \prod_{(a,b)\in I(g^{-1})}X_{a\,b}\right.\\
&\left. \hspace{32pt} \cdot \prod_{(a,b)\in I(g_{1}^{-1})\backslash I(g^{-1})} X_{b\,a}\right) g \\
&=\left(\prod_{(a,b)\in I(g_{1}^{-1})\backslash I(g^{-1})}X_{a\,b}\, \cdot \prod_{(a,b)\in I(g_{1}^{-1})\backslash I(g^{-1})}X_{b\,a}\, \cdot \prod_{(a,b)\in I(g^{-1})}X_{a\,b}\right) g \\
&=\prod_{(a,b)\in I(g_{1}^{-1})\backslash I(g^{-1})}X_{\{a,\,b\}}\, \cdot \left(\prod_{(a,b)\in I(g^{-1})}X_{a\,b}\, g\right)=X(g_{1}, g_{2})\, g^*.
\end{align*}
Here we have used the following properties
$$\prod_{(a,b)\in I((g_{1}g_{2})^{-1})}X_{a\,b}=\prod_{(a,b)\in I((g_{1}g_{2})^{-1})\cap I(g_{1}^{-1})}X_{a\,b}\, \cdot \prod_{(a,b)\in I((g_{1}g_{2})^{-1})\backslash I(g_{1}^{-1})}X_{a\,b},$$
$$\prod_{(a,b)\in I(g_{1}^{-1})}X_{a\,b}=\prod_{(a,b)\in I(g_{1}^{-1})\backslash I((g_{1}g_{2})^{-1})}X_{a\,b}\, \cdot \prod_{(a,b)\in I(g_{1}^{-1})\cap I((g_{1}g_{2})^{-1})}X_{a\,b}$$
and the proof is finished.
\begin{Corollary}\label{mig}
\begin{equation}\label{m1} 
g_{1}^*\cdot g_{2}^*=(g_{1}g_{2})^* \hspace{15pt} if \hspace{10pt} l(g_{1}g_{2})=l(g_{1})+l(g_{2})
\end{equation}
where \, $l(g):=Card\, I(g)$ \, is the lenght of\, $g\in S_{n}$. 
\end{Corollary}
\textbf{Proof.}\hspace{5pt} It is easy to see that in the case $l(g_{1})+l(g_{2})=l(g_{1}g_{2})$ we have $X(g_{1}, g_{2})=1$, so (\ref{mid}) implies (\ref{m1}).\\

The factor $X(g_{1}, g_{2})$ takes care of the reduced number of inversions in the group product of $g_{1}, g_{2} \in S_{n}$.
\begin{Example} 
Let \, $g_{1}=132$, $g_{2}=312\in S_{3}$. Then\, $g_{1}g_{2}=213$, $l(g_{1})=1$, $l(g_{2})=2$, $l(g_{1}g_{2})=1$. Note that \, $g_{1}^{-1}=132$, $g_{2}^{-1}=231$, so\\
\indent $g_{1}^*\cdot g_{2}^*=\left( X_{2\,3}\, g_{1} \right)\cdot \left( X_{1\,3} X_{2\,3}\, g_{2} \right)=X_{2\,3} X_{1\,2} X_{3\,2}\, g_{1}g_{2}=X_{\{2,\,3\}} X_{1\,2}\, g_{1}g_{2}$.\\
On the other hand we have: \, $(g_{1}g_{2})^*=X_{1\,2}\, g_{1}g_{2}$, since $(g_{1}g_{2})^{-1}=213$.\\ 
Thus we get\, $g_{1}^*\cdot g_{2}^*=X_{\{2,\,3\}}\, (g_{1}g_{2})^*$\, and\, $X(g_{1}, g_{2})=X_{\{2,\,3\}}$.
\end{Example}
\begin{Example} 
For $g_{1}=132$, $g_{2}=231$ we have\, $g_{1}g_{2}=321$, $l(g_{1})=1$, $l(g_{2})=2$, $l(g_{1}g_{2})=3$. 
Further\, $g_{1}^{-1}=132$, $g_{2}^{-1}=312$ and $(g_{1}g_{2})^{-1}=321$, so we get:\\
\indent $g_{1}^*\cdot g_{2}^*=\left( X_{2\,3}\, g_{1} \right)\cdot \left( X_{1\,2} X_{1\,3}\, g_{2} \right)=X_{2\,3} X_{1\,3} X_{1\,2}\, g_{1}g_{2}$,\\
\indent $(g_{1}g_{2})^*=X_{1\,2} X_{1\,3} X_{2\,3}\, g_{1}g_{2}$.\\
Thus\, $g_{1}^*\cdot g_{2}^*=(g_{1}g_{2})^*$\, and\, $X(g_{1}, g_{2})=1$.\\
\end{Example}
We denote by $t_{a,b}$, $1\le a\le b\le n$ the following cyclic permutation in $S_{n}$\\
\begin{equation}\label{tab} 
t_{a,b}(k):=\left\{ \begin{array}{lc}
                  k & 1\le k\le a-1 \textrm{ \, or \, } b+1\le k\le n\\
                  b & k=a\\
                  k-1 & a+1\le k\le b
         \end{array}\right.
\end{equation}
which maps $b$ to $b-1$ to $b-2$ $\cdots$ to $a$ to $b$ and fixes all $1\le k\le a-1$ and $b+1\le k\le n$ (compare with notation of $t_{a,b}$ in \cite{4}).\\ 

\noindent Let $t_{b,a}$ denote the inverse of $t_{a,b}$. Then\\
\begin{equation}\label{tba} 
t_{b,a}(k):=\left\{ \begin{array}{lc}
                  k & 1\le k\le a-1 \textrm{ \, or \, } b+1\le k\le n\\
                  k+1 & a\le k\le b-1\\
                  a & k=b.
         \end{array}\right.
\end{equation}
Then the sets of inversions are given by
$${I(t_{a,b})=\{(a,j)\mid a+1\le j\le b\}},$$
$${I(t_{b,a})=\{(i,b)\mid a \le i \le b-1\}},$$ 
so the corresponding elements in ${\mathcal{A}(S_{n})}$ have the form
\begin{equation}\label{tabtw}
t_{a,b}^*=\left( \prod_{a \le i \le b-1} X_{i\,b}\right) t_{a,b}
\end{equation}
\begin{equation}\label{tbatw}
t_{b,a}^*=\left( \prod_{a+1 \le j \le b} X_{a\,j}\right) t_{b,a}.
\end{equation}
\begin{Remark}\label{cases}
Observe that if \, $b=a$\, then\, $t_{a,a}^*=id$, $\left( \textrm{where\, } {I(t_{a,a})=\emptyset}\right)$.\\ 
In the case\, $b=a+1$\, we have\, $t_{a,a+1}=t_{a+1,a}$ and we also denote it by \, $t_{a}(=t_{a,a+1})$, $1\le a\le n-1$ $($the transposition of adjacent letters $a$ and $a+1)$.\\       
Now it is easy to see that\, $t_{a}^*=X_{a\,a+1}\, t_{a}$, with\, ${I(t_{a})=\{(a,a+1)\}}$. 
\end{Remark}
The Theorem~\ref{tmmi} implies the following more specific properties that will be presented in the following four Corollaries.
\begin{Corollary}
For each $1\le a\le n-1$ we have
\begin{equation}\label{sqta} 
(t_{a}^*)^{2}=X_{\{a,\,a+1\}}\, id.
\end{equation}
\end{Corollary}
Here we have used that\, $t_{a}t_{a}=id$\, and\, ${X_{\{a,\,a+1\}}=X_{a\,a+1}\cdot X_{a+1\,a}}$.
\begin{Corollary}[Braid relations]\label{braid}
We have
\begin{itemize} 
\item [$(i)$] \, $t_{a}^* \cdot t_{a+1}^* \cdot t_{a}^* = t_{a+1}^* \cdot t_{a}^* \cdot t_{a+1}^*$ \hspace{7pt} for each \, $1\le a\le n-2$,
\item [$(ii)$] \, $t_{a}^* \cdot t_{b}^* = t_{b}^* \cdot t_{a}^*$ \hspace{7pt} for each \, $1\le a,b\le n-1$\, with\, $|a-b|\ge 2$.
\end{itemize}
\end{Corollary}
\begin{Corollary}\label{gtba}
For each $g\in S_{n}$, $1\le a < b\le n$ we have 
$$g^* \cdot t_{b,a}^* = \left( \prod_{a < j\le b,\, g(a) > g(j)} X_{\{g(a),\,g(j)\}} \right) (gt_{b,a})^*.$$
In the case\, $g\in S_{j}\times S_{n-j}$, $1\le j\le k\le n$\, we have
\begin{equation}\label{gtkj} 
g^* \cdot t_{k,j}^* = (gt_{k,j})^*.
\end{equation}
\end{Corollary}
Compare (\ref{gtkj}) with Corollary~\ref{mig}.
\begin{Corollary}[Commutation rules]\label{comrul}
We have
\begin{itemize} 
\item [$(i)$] \, $t_{m,k}^* \cdot t_{p,k}^* = (t_{k}^*)^{2} \cdot t_{p,k+1}^* \cdot t_{m-1,k}^*$\hspace{20pt} if \hspace{5pt} $1\le k\le m < p\le n$.
\item [$(ii)$] \, Let\, $w_{n}(=n\, n-1 \cdots 2\, 1)$ be the longest permutation in $S_{n}$. Then for every $g\in S_{n}$ we have
$$(gw_{n})^* \cdot w_{n}^* = w_{n}^* \cdot (w_{n}g)^* \left(= \prod_{a < b,\, g^{-1}(a) < g^{-1}(b)} X_{\{a,\,b\}} \right) g^*.$$
\end{itemize}
\end{Corollary}

\section{Decompositions of certain canonical elements in ${\mathcal{A}(S_{n})}$}
Here we will decompose any permutation $g$ in $S_{n}$ into cycles.\\  
\indent Observe first that any permutation $g\in S_{n}$ can be represented uniquely as $g=g_{1}t_{k_{1},1}$ with $g_{1}\in S_{1}\times S_{n-1}$ and $1\le k_{1}\le n$. Then $g(k_{1})=g_{1}(t_{k_{1},1}(k_{1}))=g_{1}(1)=1$, so $k_{1}$ should be $g^{-1}(1)$.\\
Subsequently, the permutation $g_{1}\in S_{1}\times S_{n-1}$ can be represented uniquely as $g_{1}=g_{2}t_{k_{2},2}$ with $g_{2}\in S_{1}\times S_{1}\times S_{n-2}$ and $2\le k_{2}\le n$. Then $g_{1}(k_{2})=g_{2}(t_{k_{2},2}(k_{2}))=g_{2}(2)=2$ implies $k_{2}=g_{1}^{-1}(2)$.\\ 
By repeating the above procedure for every $1\le j\le n$ we can deduce that the permutation $g_{j-1}\in S_{1}^{j-1}\times S_{n-j+1}$ can be represented uniquely as $g_{j-1}=g_{j}t_{k_{j},j}$ with $g_{j}\in S_{1}^{j}\times S_{n-j}$ and $j\le k_{j}\le n$, where $g_{j-1}(k_{j})=g_{j}(t_{k_{j},j}(k_{j}))=g_{j}(j)=j$ implies $k_{j}=g_{j-1}^{-1}(j)$. Thus we get the decomposition:
\begin{equation}\label{gdecomp} 
         g =t_{k_{n},n}\cdot t_{k_{n-1},n-1}\cdots t_{k_{j},j}\cdots t_{k_{2},2}\cdot t_{k_{1},1}\, \left(= \prod_{1\le j\le n}^{\gets} t_{k_{j},j} \right). 
\end{equation}

\begin{Example} 
By applying the decomposition $(\ref{gdecomp})$ on all permutations in $S_{3}=\{123, 132, 312, 321, 231, 213 \}$ we obtain $$123=t_{3,3}t_{2,2}t_{1,1}, \hspace{10pt} 132=t_{3,3}t_{3,2}t_{1,1}, \hspace{10pt} 312=t_{3,3}t_{3,2}t_{2,1},$$
$$321=t_{3,3}t_{3,2}t_{3,1}, \hspace{10pt} 231=t_{3,3}t_{2,2}t_{3,1}, \hspace{10pt} 213=t_{3,3}t_{2,2}t_{2,1},$$
so the corresponding elements in the algebra ${\mathcal{A}(S_{3})}$ are given by
$$123^* = t_{3,3}^*\cdot t_{2,2}^*\cdot t_{1,1}^*, \hspace{10pt} 132^* = t_{3,3}^*\cdot t_{3,2}^*\cdot t_{1,1}^*, \hspace{10pt} 312^* = t_{3,3}^*\cdot t_{3,2}^*\cdot t_{2,1}^*,$$
$$321^* = t_{3,3}^*\cdot t_{3,2}^*\cdot t_{3,1}^*, \hspace{10pt} 231^* = t_{3,3}^*\cdot t_{2,2}^*\cdot t_{3,1}^*, \hspace{10pt} 213^* = t_{3,3}^*\cdot t_{2,2}^*\cdot t_{2,1}^*.$$
The following calculation shows the general situation, which will be used later in many calculations. Assume that \, ${\alpha^*_{3}} = \sum_{g\in S_{3}} g^*$. Then we get
\begin{align*}
{\alpha^*_{3}}&=\sum_{g\in S_{3}} g^* = t_{3,3}^*\cdot t_{2,2}^*\cdot t_{1,1}^* + t_{3,3}^*\cdot t_{3,2}^*\cdot t_{1,1}^* + t_{3,3}^*\cdot t_{3,2}^*\cdot t_{2,1}^*\\
&\hspace{55pt}+t_{3,3}^*\cdot t_{3,2}^*\cdot t_{3,1}^* + t_{3,3}^*\cdot t_{2,2}^*\cdot t_{3,1}^* + t_{3,3}^*\cdot t_{2,2}^*\cdot t_{2,1}^*\\
&=\left( t_{3,3}^*\right)\cdot \left( t_{3,2}^*\cdot\left( t_{3,1}^* + t_{2,1}^* + t_{1,1}^*\right) + t_{2,2}^*\cdot\left( t_{3,1}^* + t_{2,1}^* + t_{1,1}^*\right) \right)\\
&=\left( t_{3,3}^*\right)\cdot \left( t_{3,2}^* + t_{2,2}^*\right)\cdot\left( t_{3,1}^* + t_{2,1}^* + t_{1,1}^*\right)
\end{align*}
i.e
\begin{equation}\label{alpha3} 
{\alpha^*_{3}} = {\beta^*_{1}} \cdot {\beta^*_{2}} \cdot {\beta^*_{3}},
\end{equation}
where we have used the notations
\begin{align*}
{\beta^*_{1}}&=t_{3,3}^*\, \left( = id \right);\\
{\beta^*_{2}}&=t_{3,2}^* + t_{2,2}^*\, \left( = t_{3,2}^* + id \right);\\
{\beta^*_{3}}&=t_{3,1}^* + t_{2,1}^* + t_{1,1}^*\, \left( = t_{3,1}^* + t_{2,1}^* + id \right).
\end{align*} 
Therefore, we can conclude that the element ${\alpha^*_{3}}\in {\mathcal{A}(S_{3})}$ given by ${\alpha^*_{3}} = \sum_{g\in S_{3}} g^*$ can be written in the product form $(\ref{alpha3})$.
\end{Example}
In the next theorem we will prove that the element ${\alpha^*_{n}}\in {\mathcal{A}(S_{n})}$ given by ${\alpha^*_{n}} = \sum_{g\in S_{n}} g^*$, $n\ge 1$ can be decomposed into the product of simpler elements of the algebra ${\mathcal{A}(S_{n})}$ which we denote by ${\beta^*_{n-k+1}}$ for each $1\le k\le n$.
\begin{Definition}\label{beta}
For every $1\le k\le n$ we define
\begin{equation}
{\beta^*_{n-k+1}}:= t_{n,k}^* + t_{n-1,k}^* + \cdots + t_{k+1,k}^* + t_{k,k}^*\, \left(= \sum_{k\le m\le n}^{\gets} t_{m,k}^*\right).
\end{equation}
\end{Definition}
\begin{Remark}
Now it is easy to see that
\begin{align*} 
{\beta^*_{n}}&:= t_{n,1}^* + t_{n-1,1}^* + \cdots + t_{2,1}^* + t_{1,1}^* \hspace{10pt} (\textrm{if } \hspace{5pt} k = 1),\\
{\beta^*_{n-1}}&:= t_{n,2}^* + t_{n-1,2}^* + \cdots + t_{3,2}^* + t_{2,2}^* \hspace{10pt} (\textrm{if } \hspace{5pt} k = 2),\\
\vdots\\
{\beta^*_{3}}&:= t_{n,n-2}^* + t_{n-1,n-2}^* + t_{n-2,n-2}^* \hspace{10pt} (\textrm{if } \hspace{5pt} k = n-2),\\
{\beta^*_{2}}&:= t_{n,n-1}^* + t_{n-1,n-1}^* \hspace{10pt} (\textrm{if } \hspace{5pt} k = n-1),\\
{\beta^*_{1}}&:= t_{n,n}^* ( = id) \hspace{10pt} (\textrm{if } \hspace{5pt} k = n).
\end{align*}
\end{Remark}
\begin{Theorem}\label{alpha} 
Let ${\alpha^*_{n}}$ be the following canonical element in ${\mathcal{A}(S_{n})}:$
\begin{equation}\label{alphan} 
         {\alpha^*_{n}}=\sum_{g\in S_{n}} g^*.
\end{equation}
Then ${\alpha^*_{n}}$ has the following factorization 
\begin{equation}\label{alphafact} 
         {\alpha^*_{n}} = {\beta^*_{1}} \cdot {\beta^*_{2}} \cdots {\beta^*_{n}}\, \left( =\prod_{1\le k\le n}^{\gets} {\beta^*_{n-k+1}}\right).
\end{equation}
\end{Theorem}
\textbf{Proof.}\hspace{5pt} By considering decomposition (\ref{gdecomp}) of $g\in S_{n}$ and the property (\ref{gtkj}) we can write:
\begin{align*}
{\alpha^*_{n}}&=\sum_{g\in S_{n}} g^* = \sum_{\substack{g_{1}\in S_{1}\times S_{n-1}\\ 1\le k_{1}\le n}} (g_{1}t_{k_{1},1})^* = \sum_{\substack{g_{1}\in S_{1}\times S_{n-1}\\ 1\le k_{1}\le n}} g_{1}^* t_{k_{1},1}^*\\
&=\left( \sum_{g_{1}\in S_{1}\times S_{n-1}} g_{1}^*\right)\cdot \left( \sum_{1\le k_{1}\le n} t_{k_{1},1} \right) = \left( \sum_{\substack{g_{2}\in S_{1}^{2}\times S_{n-2}\\ 2\le k_{2}\le n}} (g_{2}t_{k_{2},2})^*\right)\cdot \left( \sum_{1\le k_{1}\le n} t_{k_{1},1}^* \right)\\
&=\left( \sum_{g_{1}\in S_{1}\times S_{n-1}} g_{1}^*\right)\cdot \left( \sum_{2\le k_{2}\le n} t_{k_{2},2}^* \right)\cdot \left( \sum_{1\le k_{1}\le n} t_{k_{1},1}^* \right)=\cdots=\\
&=\left( t_{k_{n},n}^*\right) \cdot \left( \sum_{n-1\le k_{n-1}\le n} t_{k_{n-1},n-1}^* \right)\cdots \left( \sum_{2\le k_{2}\le n} t_{k_{2},2}^* \right)\cdot \left( \sum_{1\le k_{1}\le n} t_{k_{1},1}^* \right)\\
&=\prod_{1\le k\le n}^{\gets} {\beta^*_{n-k+1}}
\end{align*} 
and the proof is finished.\\

Let us introduce some new elements in the algebra ${\mathcal{A}(S_{n})}$ by which we will reduce ${\beta^*_{n-k+1}}$, $1\le k\le n-1$. The motivation is to show that the element ${\alpha^*_{n}}\in {\mathcal{A}(S_{n})}$ can be expressed in turn as products of yet simpler elements of the algebra ${\mathcal{A}(S_{n})}$.  
\begin{Definition}\label{gade}
For every\, $1\le k\le n-1$ we define the following elements in the algebra ${\mathcal{A}(S_{n})}$
\begin{align*}
{\gamma^*_{n-k+1}}&:=\left( id-t_{n,k}^*\right)\cdot \left( id-t_{n-1,k}^*\right)\cdots
\left( id-t_{k+1,k}^*\right)= \prod_{k+1\le m\le n}^{\gets} \left(id-t_{m,k}^*\right),\\
{\delta^*_{n-k+1}}&:=\left( id-(t_{k}^*)^{2} \, t_{n,k+1}^*\right)\cdot \left( id-(t_{k}^*)^{2} \, t_{n-1,k+1}^*\right)\cdots
\left( id-(t_{k}^*)^{2} \, t_{k+1,k+1}^*\right)\\
&=\prod_{k+1\le m\le n}^{\gets} \left(id-(t_{k}^*)^{2} \, t_{m,k+1}^*\right)
\end{align*}
where $(t_{k}^*)^{2}$ is given by $(\ref{sqta})$ and\, $t_{k+1,k+1}^* = id$.
\end{Definition}
\begin{Proposition}\label{factorbeta} 
For every $1\le k\le n-1$ we have the following factorization
$${\beta^*_{n-k+1}} = {\delta^*_{n-k+1}} \cdot \left( {\gamma^*_{n-k+1}}\right)^{-1}.$$
\end{Proposition}
\textbf{Proof.}\hspace{5pt} Let 
$${\beta^*_{n-k+1,p}}:= \sum_{k\le m\le p}^{\gets} t_{m,k}^*$$ 
for every \, $k\le p\le n.$ Then we obtain:
\begin{align*}
{\beta^*_{n-k+1,p}}\cdot \left( id-t_{p,k}^* \right)&=\sum_{k\le m\le p}^{\gets} t_{m,k}^* - \sum_{k\le m\le p}^{\gets} t_{m,k}^* \, t_{p,k}^*\\ 
&=t_{p,k}^* + \sum_{k\le m\le p-1}^{\gets} t_{m,k}^* - \sum_{k+1\le m\le p}^{\gets} t_{m,k}^* \, t_{p,k}^* - t_{k,k}^* \, t_{p,k}^*\\
&=\sum_{k\le m\le p-1}^{\gets} t_{m,k}^* - \sum_{k+1\le m\le p}^{\gets} t_{m,k}^* \, t_{p,k}^*\\
&=\sum_{k\le m\le p-1}^{\gets} t_{m,k}^* - \sum_{k+1\le m\le p}^{\gets} (t_{k}^*)^{2} \, t_{p,k+1}^* \, t_{m-1,k}^*\\
&=\sum_{k\le m\le p-1}^{\gets} t_{m,k}^* - \sum_{k\le m\le p-1}^{\gets} (t_{k}^*)^{2} \, t_{p,k+1}^* \, t_{m,k}^*\\
&=\sum_{k\le m\le p-1}^{\gets} \left(id - (t_{k}^*)^{2} \, t_{p,k+1}^*\right) \cdot t_{m,k}^*\\
&=\left(id - (t_{k}^*)^{2} \, t_{p,k+1}^*\right) \cdot {\beta^*_{n-k+1,p-1}}
\end{align*} 
i.e
\begin{equation}\label{beta-p}
{\beta^*_{n-k+1,p}}\cdot \left( id-t_{p,k}^* \right)=\left(id - (t_{k}^*)^{2} \, t_{p,k+1}^*\right)\cdot {\beta^*_{n-k+1,p-1}}
\end{equation}
for every\, $k\le p\le n$. Note that \, ${\beta^*_{n-k+1,k}} = id$ \, and \, ${\beta^*_{n-k+1,n}} = {\beta^*_{n-k+1}}$, so for $p = n$ the identity (\ref{beta-p}) is given by
\begin{equation}\label{beta-n}
{\beta^*_{n-k+1}}\cdot \left( id-t_{n,k}^* \right)=\left(id - (t_{k}^*)^{2} \, t_{n,k+1}^*\right)\cdot {\beta^*_{n-k+1,n-1}}.
\end{equation}
By multiplying (\ref{beta-n}) from right to left with $\left(id - t_{n-1,k}^*\right)\cdots\left(id - t_{k+2,k}^*\right)\cdot\left(id - t_{k+1,k}^*\right)$ and by using above identities (\ref{beta-p}) for all $k\le p\le n-1$ it is easy to check that
\begin{align*}
&{\beta^*_{n-k+1}}\cdot \left( id-t_{n,k}^* \right)\cdot \left(id - t_{n-1,k}^*\right)\cdots\left(id - t_{k+2,k}^*\right)\cdot\left(id - t_{k+1,k}^*\right)\\
&\hspace{15pt}=\left(id - (t_{k}^*)^{2} \, t_{n,k+1}^*\right) \cdots \left(id - (t_{k}^*)^{2} \, t_{k+2,k+1}^*\right)\cdot \left(id - (t_{k}^*)^{2} \right) 
\end{align*} 
i.e 
$${\beta^*_{n-k+1}}\cdot {\gamma^*_{n-k+1}} = {\delta^*_{n-k+1}}$$
for every $1\le k\le n-1$ whence arises the identity of the Proposition~\ref{factorbeta}.
\begin{Example}\label{alpha234}
By applying $(\ref{alphafact})$ and Proposition~$\ref{factorbeta}$ we will illustrate the factorization of ${\alpha^*_{n}}\in{\mathcal{A}(S_{n})}$ in cases $n=2,3,4$\, $($recall that ${\beta^*_{1} = id})$. 
\begin{itemize}
\item [$(i)$] In the case $n=2$ we have \, ${\alpha^*_{2}} = {\beta^*_{2}}$\, and
$${\alpha^*_{2}}=\left(id-(t_{1}^*)^{2}\right)\cdot \left( id-t_{2,1}^*\right)^{-1}.$$ 
\item [$(ii)$] For $n=3$ \, we have \, ${\alpha^*_{3}} = {\beta^*_{2}}\cdot {\beta^*_{3}}$, where
\begin{align*}
{\beta^*_{2}}&=\left(id-(t_{2}^*)^{2}\right)\cdot \left( id-t_{3,2}^*\right)^{-1},\\
{\beta^*_{3}}&=\left(id-(t_{1}^*)^{2}\cdot t_{3,2}^*\right)\cdot \left(id-(t_{1}^*)^{2}\right)\cdot \left( id-t_{2,1}^*\right)^{-1}\cdot \left( id-t_{3,1}^*\right)^{-1}.
\end{align*}
\item [$(iii)$] For $n=4$ \, we have \, ${\alpha^*_{4}} = {\beta^*_{2}}\cdot {\beta^*_{3}}\cdot {\beta^*_{4}}$, where
\begin{align*}
{\beta^*_{2}}&=\left(id-(t_{3}^*)^{2}\right)\cdot \left( id-t_{4,3}^*\right)^{-1},\\
{\beta_{3}}&=\left(id-(t_{2}^*)^{2}\cdot t_{4,3}^*\right)\cdot  \left(id-(t_{2}^*)^{2}\right)\cdot \left( id-t_{3,2}^*\right)^{-1}\cdot \left( id-t_{4,2}^*\right)^{-1},\\ 
{\beta^*_{4}}&=\left(id-(t_{1}^*)^{2}\cdot t_{4,2}^*\right)\cdot \left(id-(t_{1}^*)^{2}\cdot t_{3,2}^*\right)\cdot \left(id-(t_{1}^*)^{2}\right)\cdot \left( id-t_{2,1}^*\right)^{-1}\\
&\hspace{15pt} \cdot \left( id-t_{3,1}^*\right)^{-1}\cdot \left( id-t_{4,1}^*\right)^{-1}.
\end{align*} 
\end {itemize} 
\end{Example}
\begin{Lemma}\label{lemat} We have
\begin{itemize}
\item [$(i)$] \hspace{3pt} $t_{b,a}^* = t_{1,n}\cdot t_{b+1,a+1}^*\cdot t_{n,1}$, \hspace{15pt} $1\le a\le b\le n$,
\item [$(ii)$] \hspace{3pt} $X_{\{a,\,a+1\}}\, id = t_{1,n}\cdot X_{\{a+1,\,a+2\}}\cdot t_{n,1}$, \hspace{15pt} $1\le a\le n-1$.
\end{itemize}
\end{Lemma}
\textbf{Proof.}
\begin{itemize}
\item [$(i)$] By (\ref{tab}), (\ref{tba}) and (\ref{tbatw}) we get
\begin{align*}
t_{1,n}\cdot t_{b+1,a+1}^*\cdot t_{n,1}&=t_{1,n}\cdot \left( \prod_{a+2 \le j \le b+1} X_{a+1\,j}\right) t_{b+1,a+1}\cdot t_{n,1}\\
&=\left( \prod_{a+1 \le j \le b} X_{a\,j}\right) t_{1,n}\cdot t_{b+1,a+1}\cdot t_{n,1} = t_{b,a}^*.
\end{align*}
Here we have used \, $t_{b,a} = t_{1,n} t_{b+1,a+1}t_{n,1}$\, (recall that\, $t_{1,n} = t_{n,1}^{-1}$).
\item [$(ii)$] Directly from the definition of $t_{1,n}:$
$$t_{1,n}\cdot X_{\{a+1,\,a+2\}}\cdot t_{n,1} = X_{\{a,\,a+1\}}\cdot t_{1,n}\cdot t_{n,1} = X_{\{a,\,a+1\}}\, id.$$
\end{itemize}
(This is equivalent to \, $(t_{a}^*)^{2} = t_{1,n}\cdot (t_{a+1}^*)^{2}\cdot t_{n,1}$).
\begin{Remark}
The elements\, ${\delta^*_{n-k+1}}\in {\mathcal{A}(S_{n})}$, $1\le k\le n-1$\, from Definition~$\ref{gade}$ can be rewritten as: 
\begin{align*}
{\delta^*_{n-k+1}}&=\left( id-X_{\{k,\,k+1\}}\, t_{n,k+1}^*\right)\cdot \left( id-X_{\{k,\,k+1\}}\, t_{n-1,k+1}^*\right)\cdots \\
&\hspace{25pt} \left( id-X_{\{k,\,k+1\}}\, t_{k+2,k+1}^*\right)\cdot \left( id-X_{\{k,\,k+1\}}\, t_{k+1,k+1}^*\right) 
\end{align*}
or shorter
\begin{equation}\label{deltak}
{\delta^*_{n-k+1}} = \prod_{k+1\le m\le n}^{\gets} \left(id-X_{\{k,\,k+1\}}\, t_{m,k+1}^*\right).
\end{equation}
Our next goal is to give a formula for the inverse of ${\alpha^*_{n}}$. In order to do this we first need to determine the inverse of \, ${\delta^*_{n-k+1}}$\, for all\, $1\le k\le n-1$, because 
$$\left({\alpha^*_{n}}\right)^{-1} = {\gamma^*_{n}} \cdot \left( {\delta^*_{n}}\right)^{-1}\cdot \, {\gamma^*_{n-1}} \cdot \left( {\delta^*_{n-1}}\right)^{-1}\cdots \, {\gamma^*_{2}} \cdot \left( {\delta^*_{2}}\right)^{-1}.$$
\end{Remark}

\noindent Let us introduce a more accurate label
\begin{equation}\label{deltank}
{\delta^*_{n-k+1,n}}:= {\delta^*_{n-k+1}}
\end{equation}
where ${\delta^*_{n-k+1}}$ is given by (\ref{deltak}).\\

\noindent Let us denote by
        $$Des(\sigma):=\{1\le i\le n-1 \mid \sigma(i) > \sigma(i+1) \}$$ 
the descent set of a permutation $\sigma\in S_{n}$.\\ 
Let \, $des(\sigma)=Card(Des(\sigma))$\, be the number of descents of $\sigma$.\\ 
Note that for $g\in S_{1}^k\times S_{n-k}$
$$Des(g)=\{k+1\le i\le n-1 \mid g(i) > g(i+1) \}.$$
\begin{Proposition}\label{invdel}
The inverse of \, ${\delta^*_{n-k+1,n}}$, $1\le k\le n-1$\, is given by the formula
\begin{equation}\label{invdeln}
\left({\delta^*_{n-k+1,n}}\right)^{-1} = \left( \Delta_{n-k+1,n} \right)^{-1}\cdot \left({\varepsilon^*_{n-k+1,n}}\right)
\end{equation}  
where
$$\Delta_{n-k+1,n} := \left( id-X_{\{k,\,k+1\}} \right)\cdot \left( id-X_{\{k,\,k+1\,k+2\}} \right)\cdots \left( id-X_{\{k,\,k+1,\,\dots ,\,n\}} \right),$$ 
$${\varepsilon^*_{n-k+1,n}} := \sum_{g\in S_{1}^k\times S_{n-k}} \omega_{n-k+1,n}(g)\, g^*$$ 
and
$$\omega_{n-k+1,n}(g) := \prod_{i\in Des(g^{-1})} X_{\{k,\,k+1,\,\dots ,\,i\}}.$$
\end{Proposition}
\textbf{Proof.}\hspace{5pt} By (\ref{deltak}) and (\ref{deltank}) we have
$${\delta^*_{n-k+1,n}} = \left(id-X_{\{k,\,k+1\}}\, t_{n,k+1}^*\right) \cdot \prod_{k+1\le m\le n-1}^{\gets} \left(id-X_{\{k,\,k+1\}}\, t_{m,k+1}^*\right)$$
or shortly
\begin{equation}\label{tdeltank}
{\delta^*_{n-k+1,n}} = \left(id-X_{\{k,\,k+1\}}\, t_{n,k+1}\right)\cdot {\delta^*_{n-k+1,n-1}}
\end{equation} 
where
\begin{align*}
{\delta^*_{n-k+1,n-1}} &= \prod_{k+1\le m\le n-1}^{\gets} \left(id-X_{\{k,\,k+1\}}\, t_{m,k+1}^*\right)\\
&= t_{1,n}\cdot \left( \prod_{k+1\le m\le n-1}^{\gets} \left(id - X_{\{k+1,\,k+2\}}\, t_{m+1,k+2}^*\right)\right)\cdot t_{n,1}\\
&= t_{1,n}\cdot \left( \prod_{k+2\le m\le n}^{\gets} \left(id - X_{\{k+1,\,k+2\}}\, t_{m,k+2}^*\right)\right)\cdot t_{n,1} = t_{1,n}\cdot {\delta^*_{n-k,n}} \cdot t_{n,1}.
\end{align*}
Here we have used property $(ii)$ of the Lemma~\ref{lemat}. Thus we obtain
$${\delta^*_{n-k+1,n}} = \left(id - X_{\{k,\,k+1\}}\, t_{n,k+1}^*\right)\cdot t_{1,n}\cdot {\delta^*_{n-k,n}} \cdot t_{n,1}$$
i.e the identity
$$\left({\delta^*_{n-k+1,n}}\right)^{-1} \cdot \left(id - X_{\{k,\,k+1\}}\, t_{n,k+1}^*\right) = t_{1,n}\cdot \left({\delta^*_{n-k,n}}\right)^{-1} \cdot t_{n,1}$$
which takes the form:
$$\left( \Delta_{n-k+1,n} \right)^{-1}\cdot {\varepsilon^*_{n-k+1,n}}\cdot \left(id - X_{\{k,\,k+1\}}\, t_{n,k+1}^*\right) = t_{1,n}\cdot \left( \Delta_{n-k,n} \right)^{-1}\cdot {\varepsilon^*_{n-k,n}} \cdot t_{n,1}$$
or
\begin{equation}\label{epsin}
{\varepsilon^*_{n-k+1,n}}\cdot \left(id - X_{\{k,\,k+1\}}\, t_{n,k+1}^*\right)=\left(id - X_{\{k,\,k+1,\,\dots ,\,n\}}\right) \cdot {\varepsilon^*_{n-k+1,n-1}}
\end{equation} 
where
\begin{align*}
{\varepsilon^*_{n-k+1,n-1}} &= t_{1,n}\cdot {\varepsilon^*_{n-k,n}}\cdot t_{n,1},\\
id - X_{\{k,\,k+1,\,\dots ,\,n\}} &= \Delta_{n-k+1,n}\cdot t_{1,n}\cdot \left(\Delta_{n-k,n}\right)^{-1}\cdot t_{n,1}. 
\end{align*}
To prove the formula (\ref{invdeln}) (by induction), it suffices to prove the identity (\ref{epsin}). 
Notice that (\ref{epsin}) is equivalent to 
$${\varepsilon^*_{n-k+1,n}} = \left(id - X_{\{k,\,k+1,\,\dots ,\,n\}}\right) \cdot {\varepsilon^*_{n-k+1,n-1}}\cdot \left(id - X_{\{k,\,k+1\}}\, t_{n,k+1}^*\right)^{-1}.$$
We first calculate
\begin{align*}
{\varepsilon^*_{n-k+1,n}}\cdot X_{\{k,\,k+1\}}\, t_{n,k+1}^* &= \sum_{\sigma\in S_{1}^k\times S_{n-k}} \omega_{n-k+1,n}(\sigma)\, {\sigma}^*\cdot X_{\{k,\,k+1\}}\, t_{n,k+1}^*\\
&= \sum_{\sigma\in S_{1}^k\times S_{n-k}} \omega_{n-k+1,n}(\sigma)\cdot X_{\{k,\,\sigma(k+1)\}}\, {\sigma}^*\cdot t_{n,k+1}^*\\  
&= \sum_{\sigma\in S_{1}^k\times S_{n-k}} \omega_{n-k+1,n}(\sigma)\cdot X_{\{k,\,\sigma(k+1)\}}\\
&\hspace{20pt} \cdot \prod_{k+1\le j < \sigma(k+1)} X_{\{j,\,\sigma(k+1)\}}\, ({\sigma}\, t_{n,k+1})^*\\  
&= \sum_{\sigma\in S_{1}^k\times S_{n-k}} \omega_{n-k+1,n}(\sigma)\cdot \prod_{k\le j < \sigma(k+1)} X_{\{j,\,\sigma(k+1)\}}\, ({\sigma}\, t_{n,k+1})^*\\  
&= \sum_{g\in S_{1}^k\times S_{n-k}} \omega_{n-k+1,n}(g\, t_{n,k+1}^{-1})\cdot \prod_{k\le j < g(n)} X_{\{j,\,g(n)\}}\, g^*  
\end{align*}
where we used that \, $g=\sigma\, t_{n,k+1}$\, implies \, $\sigma = g\, t_{n,k+1}^{-1}$,\, so \, $\sigma(k+1) = g(n)$.\\
On the other hand, by the formula
\begin{equation}\label{Des} 
Des(t_{n,k+1}\, g^{-1})= \left(  Des(g^{-1})\backslash \left\{ g(n)\right\} \right) \cup \left\{ g(n)-1\right\} \hspace{18pt} \textrm{if} \hspace{10pt} g(n) \le n
\end{equation} 
where \hspace{15pt} 
\begin{equation*}
Des(g^{-1})\backslash \left\{ g(n)\right\}=Des(g^{-1}) \hspace{15pt} \textrm{when} \hspace{10pt} g(n)=n, \hspace{7pt} g(n)\notin Des(g^{-1})
\end{equation*} 
we obtain
\begin{align*}
\omega_{n-k+1,n}(g\, t_{n,k+1}^{-1})&\cdot \prod_{k\le j < g(n)} X_{\{j,\,g(n)\}}\\
&= \prod_{i\in Des(t_{n,k+1}\, g^{-1})} X_{\{k,\,k+1,\,\dots ,\,i\}}\cdot \prod_{k\le j < g(n)} X_{\{j,\,g(n)\}}\\
&=\sum_{\substack{i\in Des(t_{n,k+1}\, g^{-1})\\i\ne g(n)-1}} X_{\{k,\,k+1,\,\dots ,\,i\}}\cdot \underbrace{X_{\{k,\,k+1,\,\dots ,\,g(n)-1\}}}_{if \hspace{5pt} i=g(n)-1}\cdot \prod_{k\le j < g(n)} X_{\{j,\,g(n)\}}\\   
&=\sum_{\substack{i\in Des(t_{n,k+1}\, g^{-1})\\i\ne g(n)-1}} X_{\{k,\,k+1,\,\dots ,\,i\}}\cdot X_{\{k,\,k+1,\,\dots ,\,g(n)\}}\\
&=\left\{ \begin{array}{ll}
                  \omega_{n-k+1,n}(g) & \textrm{if \, } g(n)< n\\
                  X_{\{k,\,k+1,\,\dots ,\,n\}}\cdot \omega_{n-k+1,n}(g) & \textrm{if \, } g(n)= n.
         \end{array}\right.
\end{align*}
Therefore
\begin{align*}
{\varepsilon^*_{n-k+1,n}}&\cdot X_{\{k,\,k+1\}}\, t_{n,k+1}^*\\
&=\left\{ \begin{array}{ll}
                  \sum_{g\in S_{1}^k\times S_{n-k}} \omega_{n-k+1,n}(g)\, g^* & \textrm{if \, $g(n)< n$}\\
                  \sum_{g\in S_{1}^k\times S_{n-k}} X_{\{k,\,k+1,\,\dots ,\,n\}}\cdot\omega_{n-k+1,n}(g)\, g^* & \textrm{if \, $g(n)= n$}
         \end{array}\right.\\
&=\left\{ \begin{array}{ll}
                  {\varepsilon^*_{n-k+1,n}} & \textrm{if \, } g(n)< n\\
                  X_{\{k,\,k+1,\,\dots ,\,n\}}\cdot {\varepsilon^*_{n-k+1,n}} & \textrm{if \, } g(n)= n.
         \end{array}\right.
\end{align*}
Finally, we get
\begin{align*}
{\varepsilon^*_{n-k+1,n}}&\cdot \left(id - X_{\{k,\,k+1\}}\, t_{n,k+1}^*\right)=
{\varepsilon^*_{n-k+1,n}} - {\varepsilon^*_{n-k+1,n}}\cdot X_{\{k,\,k+1\}}\, t_{n,k+1}^*\\
&=\sum_{\substack{g\in S_{1}^k\times S_{n-k}\\ g(n)<n}} \omega_{n-k+1,n}(g)\, g^* + \sum_{\substack{g\in S_{1}^k\times S_{n-k}\\ g(n)=n}} \omega_{n-k+1,n}(g)\, g^*\\
&\hspace{15pt} - \sum_{\substack{g\in S_{1}^k\times S_{n-k}\\ g(n)<n}} \omega_{n-k+1,n}(g)\, g^* - \sum_{\substack{g\in S_{1}^k\times S_{n-k}\\ g(n)=n}} X_{\{k,\,k+1,\,\dots ,\,n\}}\cdot \omega_{n-k+1,n}(g)\, g^*\\
&=\left(id - X_{\{k,\,k+1,\,\dots ,\,n\}}\right) \cdot \sum_{\substack{g\in S_{1}^k\times S_{n-k}\\ g(n)=n}} \omega_{n-k+1,n}(g)\, g^*\\
&=\left(id - X_{\{k,\,k+1,\,\dots ,\,n\}}\right) \cdot {\varepsilon^*_{n-k+1,n-1}}
\end{align*} 
where we have used that
\begin{align*}
\sum_{\substack{g'\in S_{1}^k\times S_{n-k}\\ g'(n)=n}} \omega_{n-k+1,n}(g')\, (g')^*&=\sum_{\substack{g'\in S_{1}^k\times S_{n-k}\\ g'(n)=n}} t_{1,n}\cdot \left( t_{n,1}\cdot \omega_{n-k+1,n}(g')\, (g')^*\cdot t_{1,n}\right) \cdot t_{n,1}\\
&=t_{1,n}\cdot \left( \sum_{g\in S_{1}^{k+1}\times S_{n-k-1}} \omega_{n-k,n}(g)\, g^* \right)\cdot t_{n,1}\\
&=t_{1,n}\cdot {\varepsilon^*_{n-k,n}}\cdot t_{n,1}={\varepsilon^*_{n-k+1,n-1}}.
\end{align*}
This prove (\ref{epsin}) and the proof of the Proposition~\ref{invdel} is now completed.\\

The matrix factorizations from the right given in $\cite{4}$ and $\cite{5}$ one can replace by twisted algebra factorizations (from the right). But here we have presented the factorizations from the left because they are more suitable for computing constants in multiparametric algebra of noncommuting polynomials (this will be treated in a forthcoming paper).\\

\end{document}